\documentclass[journal,draftcls,onecolumn,12pt,twoside]{IEEEtran}

\usepackage[T1]{fontenc}

\usepackage{ifpdf}
\ifpdf
\else
\fi

\usepackage{cite}

\ifCLASSINFOpdf

\else

\usepackage[dvips]{graphicx}

\fi

\usepackage{amsmath}

\usepackage[cmintegrals]{newtxmath}

\makeatletter
\def\endthebibliography{%
	\def\@noitemerr{\@latex@warning{Empty `thebibliography' environment}}%
	\endlist
}
\usepackage{comment}
\usepackage{graphicx}
\usepackage{multirow}
\usepackage{hyperref} 
\usepackage{amsmath}
\usepackage{url}
\usepackage{caption}
\usepackage{epstopdf}
\usepackage{soul}
\usepackage{color}
\usepackage{cite}
\usepackage{enumerate}
\usepackage{xcolor}
\usepackage[utf8]{inputenc}
\usepackage[english]{babel}
\usepackage[section]{placeins}
\usepackage[linesnumbered,ruled,vlined]{algorithm2e}
\usepackage[font=scriptsize]{caption}
\hypersetup{linkbordercolor=red}

\newcommand{\Rn}{\mathbb{R}}
\newcommand{\Cn}{\mathbb{C}}

\usepackage{stackengine}
\usepackage{amsmath}

\def \A {{\mathbf A}}

\def \G {{\mathbf G}}
\def \H {{\mathbf H}}

\def \L {{\mathbf L}}

\def \R {{\mathbf R}}

\def \b {{\mathbf b}}

\def \d {{\mathbf d}}

\def \f {{\mathbf f}}
\def \g {{\mathbf g}}

\def \n {{\mathbf n}}
\def \r {{\mathbf r}}

\def \t {{\mathbf t}}

\def \v {{\mathbf v}}
\def \w {{\mathbf w}}
\def \x {{\mathbf x}}
\def \y {{\mathbf y}}
\def \z {{\mathbf z}}

\def \xh {{\hat{\x}}}
\def \xhat {{\xh(\x_t,\y_t,\xib_t)}}
\def \xhatb {{\xh(\x_{\t},\y_{\t},\xib_{\t})}}

\def \t {{[t]}}

\def \ft {{\tilde{f}}}
\def \fh {{\hat{f}}}
\def \L {{\hat{L}}}

\usepackage{amsmath}

\newcommand{\norm}[1]{\ensuremath{\left\|#1\right\|}} 
\providecommand{\ip}[1]{\langle#1\rangle}

\def \ep {{\varepsilon}}

\def \vc {{\check{\v}}}
\def \dt {{\tilde{\del}_\t}}

\usepackage{amsmath}

\providecommand{\norm}[1]{$\left\|#1\right\|$}
\providecommand{\ip}[1]{\langle#1\rangle}

\def \xib {{\boldsymbol \xi}}
\def \Xib {{\boldsymbol \Xi}}

\def \del {{\boldsymbol \Delta}}
\def \Ga {{\boldsymbol \Gamma}}
\def \Hc {{\mathcal H}}

\def \tr {{\text{tr}}}
\def \vec {{\text{vec}}}

\def \expect {{\mathbb E}}

\def \Rn {{\mathbb{R}}}

\def \so {{\boldsymbol \theta}}

\def \O {{\mathcal{O}}}

\newcommand{\ind}{1\hspace{-1.6mm}1} 

\def \Xc {{\mathcal X}}

\def \Xc {{\mathcal X}}

\def \del {{\boldsymbol{\delta}}}
\def \ph {{\boldsymbol{\phi}}}
\def \Hc {{\mathcal{H}}}

\newcommand{\Ek}[1]{\ensuremath{\mathbb{E}_{\xib_t}\left[#1\right]}}  

\newcommand{\Ex}[1]{\ensuremath{\mathbb{E}[#1]}} 		  

\newtheorem{assumption}{}

\newtheorem{lemma}{Lemma}

\newtheorem{theorem}{Theorem}

\def \del {{\boldsymbol{\delta}}}
\def \ph {{\boldsymbol{\phi}}}
\def \Hc {{\mathcal{H}}}

\begin{document}
	
	\title{Practical Precoding via Asynchronous Stochastic Successive Convex Approximation}
	
	\title{Practical Precoding via Asynchronous Stochastic Successive Convex Approximation}
	
	\author{ Basil M. Idrees,
		Javed Akhtar, and, 
		Ketan Rajawat, \IEEEmembership{Member, IEEE}
		
		\thanks{	
			Basil M. Idrees and Ketan Rajawat are with the Dept. of Electrical Engineering, Indian Institute of Technology Kanpur, India (e-mail: \texttt{\{basilmi, ketan\} @iitk.ac.in}). Javed Akhtar is with the Radisys India Pvt. Ltd., Bengaluru (e-mail: \texttt{javeda2309@gmail.com})}
	}
	\maketitle
	
	\begin{abstract}
		
		We consider stochastic optimization of a smooth non-convex loss function with a convex non-smooth regularizer. In the online setting, where a single sample of the stochastic gradient of the loss is available at every iteration, the problem can be solved using the proximal stochastic gradient descent (SGD) algorithm and its variants. However in many problems, especially those arising in communications and signal processing, information beyond the stochastic gradient may be available thanks to the structure of the loss function. Such extra-gradient information is not used by SGD, but has been shown to be useful, for instance in the context of stochastic expectation-maximization, stochastic majorization-minimization, and stochastic successive convex approximation (SCA) approaches. By constructing a stochastic strongly convex surrogates of the loss function at every iteration, the stochastic SCA algorithms can exploit the structural properties of the loss function and achieve superior empirical performance as compared to the SGD.
		
		In this work, we take a closer look at the stochastic SCA algorithm and develop its asynchronous variant which can be used for resource allocation in wireless networks. While the stochastic SCA algorithm is known to converge asymptotically, its iteration complexity has not been well-studied, and is the focus of the current work. The insights obtained from the non-asymptotic analysis allow us to develop a more practical asynchronous variant of the stochastic SCA algorithm which allows the use of surrogates calculated in earlier iterations. We characterize precise bound on the maximum delay the algorithm can tolerate, while still achieving the same convergence rate. We apply the algorithm to the problem of linear precoding in wireless sensor networks, where it can be implemented at low complexity but is shown to perform well in practice.
		
	\end{abstract}

	\section{Introduction}\label{intro}
	This work considers the stochastic optimization problem
	\begin{align}
	\x^\star = \arg\min_{\x\in\Xc} &~U(\x) :=  F(\x) + h(\x)     \tag{$\mathcal{P}$}  \label{opt-P4}
	\end{align}
	where $F(\x) := \mathbb{E}_{\xib}[f(\x,\xib)]$ is a smooth and possibly non-convex function, $\Xc \subset \mathbb{R}^n$ is a proper closed convex set, and  $h$ is convex but possibly non-smooth function. The expectation is with respect to the random variable $\xib \in \Rn^d$, whose distribution is not known a priori. Instead, an \emph{online} or streaming setting is considered wherein the independent identically distributed (i.i.d.) samples $\xib_t$ are observed in a sequential fashion.
	
	The stochastic optimization problem in \eqref{opt-P4} has been well studied in the context of several signal processing, communications, and machine learning applications \cite{mokhtari2017large}. Since the problems arising in these domains are often high-dimensional, thus most efficient algorithms for solving \eqref{opt-P4} only require access to a stochastic first-order oracle (SFO) that provides $\nabla f(\x_t,\xib_t)$ at every $t$. The state-of-the-art approach to solving \eqref{opt-P4} is through the use of the proximal stochastic gradient descent (prox-SGD) algorithm with updates \cite{ghadimi2016mini}
	\begin{align}\label{prox}
	\x_{t+1} =  \arg\min_{\x \in \Xc} h(\x) + \check{f}(\x,\x_t,\xib_t)
	\end{align}
	where $\check{f}(\x,\x_t): = f(\x_t) + \ip{\nabla f(\x_t,\xib_t),\x} + \frac{1}{2\eta}\norm{\x-\x_t}^2$ and  $\eta$ is the step-size parameter. The prox-SGD is particularly attractive when the per-update optimization problem \eqref{prox} can be solved efficiently. It is known that in order to reach an $\epsilon$-stationary solution to \eqref{opt-P4}, the prox-SGD incurs an SFO complexity of at least $\O(\epsilon^{-2})$ \cite{ghadimi2016mini}.
	
	The SFO assumption is quite general and works even for problems with a black-box access to function values or gradients. For many problems however, operating under such a general rubric may turn out to be too restrictive, as the performance bounds may not reveal the importance of using any information beyond the stochastic gradients. At times, the objective function may be entirely known in closed-form and may possess structure that can be exploited to obtain algorithms that are faster than SGD. For instance, Bayesian learning problems are often non-convex but have a specific block-separable structure that makes them amenable to expectation-maximization (EM) \cite{hoffman2013stochastic} class of algorithms. As a result, EM and its stochastic variants have been the methods of choice for Bayesian learning problems since several decades. Remarkably, while there are many other algorithms that also exploit the structure of the objective function, their oracle complexity analysis has not attracted much attention.
	
	In this work, we develop algorithms utilizing the stochastic convex approximation oracle (SCO), wherein at each $t$, a strongly convex surrogate function $\hat{f}(\x,\x_t,\xib_t)$, that satisfies the tangent condition $\nabla \hat{f}(\x_t,\x_t,\xib_t) = \nabla f(\x_t,\xib_t)$, is revealed. Note that if $f$ is smooth, such a surrogate always exists, and $\check{f}(\x,\x_t,\xib_t)$ used in \eqref{prox} is an example. In other words, the surrogate functions revealed by the oracle may not necessarily contain any information beyond the first order stochastic gradients, suggesting that a lower bound on the SFO complexity for \eqref{opt-P4} is also a lower bound on the SCO complexity for \eqref{opt-P4}. Nevertheless, more general surrogate functions have been widely used in the context of successive convex approximation (SCA) algorithms; see \cite{scutari2016parallel, scutari2016paralel} for examples. The stochastic SCA and its variants have likewise been proposed in \cite{Aliu23,liu2018stochastic,liu2018online,liu2019two} but only asymptotic analysis is carried out. On the other hand, the convergence rate obtained in \cite{koppel2018parallel} is suboptimal.  
	
	The main contribution of this paper is the development of SCO complexity bounds for a class of algorithms for solving \eqref{opt-P4}. As in SGD, the non-asymptotic analysis is important as it reveals the dependence of the algorithm performance on various problem and algorithm parameters. We also study the performance of the proposed algorithm under delayed updates and characterize an explicit bound on the maximum tolerable delay. The proposed asynchronous stochastic SCA algorithm can therefore be used in computationally constrained settings, where the surrogate function minimization subproblem cannot be solved instantaneously, and therefore the updates must be carried out at a later iteration. 
	
	As a second contribution, we demonstrate the usefulness of the proposed algorithm to resource allocation problems in wireless networks. The standard approach to real-time resource allocation, at least in theory, entails observing the environmental state and formulating an optimization problem that yields the appropriate resource variables. For instance, in multi-input multi-output (MIMO) wireless networks, one would observe the channel gains at every coherence interval and design the optimal precoding matrix that minimizes the mean-square error or bit-error rate. However, from the implementation stand-point, these optimization problems are often too complicated to be solved concurrently, and designers instead opt for heuristic approaches whose performance may be far from optimal. The proposed SCA framework offers a flexible resource allocation framework wherein the allocation variables are split into static and time-varying components. When the channel variations are small, the time-varying component of the resource (such as the precoder matrix) is determined quickly by solving an approximate per-iteration optimization problem, while an estimate of the static component is learned over time. The asynchronous nature of the proposed algorithm allows the ensuing updates to be carried out at the edge devices and simplifies the implementation. 
	
	\subsection{Literature Review}
	We briefly review the literature on algorithms for stochastic non-convex optimization and their asynchronous variants. We limit our discussion to online algorithms where observations arrive in a sequential manner and only a single sample is processed at every time instant. In other words, we do not discuss algorithms relying on finite-sum structure or those using mini-batching or momentum. 
	
	\subsubsection{Stochastic Non-convex Optimization} 
	When only stochastic (sub-)gradient information is available, \eqref{opt-P4} can be solved in an online fashion using the proximal SGD algorithm \cite{ghadimi2013stochastic}. Variants of prox-SGD for weakly convex objective function and regularizers  \cite{davis2019stochastic,davis2019proximally, asi2019stochastic} as well as for non-smooth non-convex regularizers \cite{xu2019non,xu2019stochastic,metel2019simple} have since been proposed. 
	
	As discussed earlier, the stochastic successive convex approximation algorithms entail constructing a stochastic surrogate of the objective function at every iteration, and the update involves solving the resulting convex optimization problem. One of the first algorithms in this class was the stochastic majorization-minimization (MM) algorithm, which was asymptotically convergent to a stationary point, but required stringent conditions on the surrogate function. Variants of the stochastic MM algorithm have been applied to several applications in signal processing; see e.g., \cite{Aliu21,chouzenoux2016convergence} and references therein. While some of these works do obtain non-asymptotic convergence rates, the analysis is often tailored to the specific problems at hand and not easy to generalize. 
	
	SCA algorithms requiring less stringent conditions on the surrogate functions were first proposed for deterministic problems in \cite{scutari2016parallel, scutari2016paralel}, and subsequently extended to various stochastic problems in \cite{yang2016parallel,liu2018stochastic,ye2019stochastic,liu2018online,liu2019two}, among others. Of these, \cite{liu2018stochastic,liu2018online,ye2019stochastic,liu2019two}  also allow convex approximation of the non-convex constraint functions, and provide asymptotic convergence guarantees to a stationary point. The SCO complexity of these algorithms for the general stochastic constrained non-convex setting have not been obtained so far. 
	
	The problem formulation in \eqref{opt-P4} and the stochastic SCA algorithm is a special case of those in \cite{liu2018stochastic}, since we assume that the constraints are deterministic and convex. The specialization is necessary as it allows us to obtain the required SCO complexity. A similar setting has also been considered in \cite{mokhtari2017large, koppel2018parallel}, where the SCO complexity of $\mathcal{O}(\epsilon^{-4})$ is derived. In contrast, the synchronous version of the proposed algorithm in this work yields an SCO complexity of $\mathcal{O}(\epsilon^{-2})$, which is at par with that of the prox-SGD algorithm. 
	
	\subsubsection{Asynchronous Algorithms}
	While asynchronous algorithms have been well-studied since several decades \cite{bertsekas1989parallel}, the corresponding rate results were only characterized in the last decade \cite{recht2011hogwild}, \cite{asgd2}. Asynchronous variants of the deterministic SCA algorithm were proposed in \cite{cannelli2019asynchronous,cannelli2017asynchronous1}. However, these results do not carry over to the stochastic setting, and consequently, there do not exist any asynchronous variants of the stochastic SCA algorithm. In this work, we proposed such an asynchronous variant, and obtain its SCO complexity. 
	
	\subsubsection{Resource Allocation in Wireless Systems}
	This work considers vector parameter estimation in wireless sensor networks (WSN), where the goal is to design the precoder matrix. The general problem is non-convex but has been well-studied in various contexts \cite{xiao2008linear,behbahani2012linear,akhtar2018distributed}. A common thread among these algorithms is the need to solve a complicated optimization problem at every coherence interval, in order to obtain the required precoding matrix. In general, vector quantization can be used to reduce the complexity of precoding, since only one of the stored precoding matrices need to be used, based on the feedback \cite{he2017codebook,fanaei2013limited}. Complementary to the quantization-based approaches, we consider the case where the channel can be written as a sum of a static large-scale fading component and the dynamic small-scale fading component. Assuming that the dynamic component of the channel is small, as is usually the case in IoT communications  \cite{savaux2019uplink,beyene2017performance}, the proposed approach constructs a static and a dynamic component of the precoder matrix in an online fashion. Such a split is commonly used in robust precoder design \cite{venkategowda2017precoding}. 
	
	More generally, two time-scale resource allocation has been considered before; see e.g. \cite{liu2019two}, but its convergence rate has not been obtained. Parallel to these, asynchronous algorithms for solving general convex  resource allocation problems in wireless systems have also been studied before \cite{bedi2018asynchronous,rajawat2011cross}. The present work is however the first to apply asynchronous stochastic SCA algorithms for solving non-convex resource allocation problems.
	
	\subsection{Notations}
	Before going further lets describe some of the notations used in this work. The bold lower case letters  are used for denoting vectors while bold upper case letters are used for denoting matrices.  $\A^\top$, $\A^H$, and  $\text{tr}(\A)$, are used for denoting the transpose, Hermitian transpose,  and trace of a matrix $\A$   respectively. Also,  $\vec(\A)$ denotes vectorized version of $\A$  which is obtained by stacking its columns into a single column vector. The expectation operation is denoted by $\expect$.  The space of all real and complex matrices of size $m \times n$ is denoted by $\mathbb{R}^{m \times n}$ and $\mathbb{C}^{m \times n}$, respectively.

	The rest of the paper is organized as follows: Section II presents the problem statement  followed by a couple of examples. Section III puts forth the proposed algorithm  and develops the required SCO bounds. Section IV applies the proposed algorithm to the precoder design problem in WSNs and discusses the numerical results. 
	
	\section{Problem Formulation}
	As discussed earlier, the problem in \eqref{opt-P4} is a special case of the general constrained problem in \cite{liu2018stochastic}. The stochastic SCA algorithm proposed in \cite{liu2018stochastic} entails carrying out the updates by solving
	\begin{align}\label{xhata0}
	\xh(\x_t,\y_t,\xib_t) &= \arg\min_{\x\in\Xc} h(\x) + 
	\ft(\x,\x_t,\y_t,\xib_t)
	\end{align}
	where $\ft$ is a convex surrogate of $f$ and $\x_t$, $\y_t$ are iterates that are  updated as
	\begin{subequations}\label{sync}
		\begin{align}
		\y_{t+1} &= (1-\rho_t)\y_t + \rho_t \nabla f(\x_t,\xib_t)\\
		\x_{t+1} &= (1-\gamma_t)\x_t + \gamma_t \xh(\x_t,\y_t,\xib_t) 
		\end{align}
	\end{subequations}
	with $\rho_t$ and $\gamma_t$ being diminishing sequences. This work will develop an asynchronous variant of \eqref{xhata0} and study its SCO complexity. 
	
	In order to motivate the need for an asynchronous algorithm, consider that for many applications, especially those arising in wireless networks, calculating the surrogate function may be costly. In such settings, considerable duration may elapse between the time when $\xib_t$ is observed and the time when $\xh(\x_t,\y_t,\xib_t)$ is calculated. Naturally, the updates in \eqref{sync}  cannot be carried out while these quantities are being calculated and the algorithm must wait. If the time required to calculate $\xh(\x_t,\y_t,\xib_t)$ exceeds a single time slot, the excess waiting time would accumulate over time, resulting in an ever increasing lag between the time $\xib_t$ is observed and the time when it is ultimately used. As a simple example, if the gradient calculation takes up time equal to 2 coherence intervals, the $t$-th iteration of the algorithm occurs at time $2t$, and the storage requirements increase in an unbounded fashion.
	
	\subsection{SCA for resource allocation}\label{envelope}
	We discuss a typical resource allocation example where the SCA updates are complicated and cannot be calculated easily. We begin with discussing a generic network resource allocation problem and subsequently provide a concrete example involving a wireless sensor network. Specifically, consider a problem where the objective function has composite structure:
	\begin{align}
	f(\x,\xib) &= \ell(\f(\x,\xib),\xib) \label{ell}\\
	f_j(\x,\xib) &= \max_{\w} f_j^o(\x,\w,\xib)  \label{fj}\\
	\text{s. t. }& \hspace{4mm} f_j^c(\x,\w,\xib) \geq 0 \hspace{3.5mm} 1\leq j \leq m \label{fj2}
	\end{align}
	where $f_j, f_j^o, f_j^c:\Rn^n \times \Rn^d \rightarrow \Rn$ and $\{f_j\}_{j=1}^m$ are the entries of the vector-valued function $\f$. Here, $\ell:\Rn^m \rightarrow \Rn$ is the outer function, while the inner functions $f_j$ are defined as maximum values of $f^o_j$ over $\w \in \Rn^s$ under the constraints in \eqref{fj2}. Observe that the optimum $\w$ in \eqref{fj} would be a function of $\x$ and $\xib$ and would therefore be different for each time instant. In wireless settings, $\w$ would include \emph{instantaneous} variables such as power allocations or entries of the beamforming/precoding/decoding matrices. On the other hand, $\x$ would include long-term allocation variables such as link rates, average power allocations, design variables, etc. The inner function, for instance, may be the rate achieved by a user subject to local power constraints, while the outer function may represent the (negative of the) system utility derived from these rates. We consider the setting where the inner subproblem is a convex optimization problem and is readily solvable. 
	
	A surrogate for $f$ can be constructed through the use of the envelope theorem. Specifically, let $(\w_j^\star,\lambda_j^\star)$ denote the primal-dual optimum of the inner subproblem \eqref{fj} for a given $\x$ and $\xib$, i.e., $(\w_j^\star,\lambda_j^\star) = \arg\min_{\lambda}\max_{\w} f_j^o(\x,\w,\xib)+\lambda f_j^c(\x,\w,\xib)$ for all $j = 1, \ldots, m$. Then, from the envelope theorem, the gradient of $f_j$ is given by 
	\begin{align}
	\nabla_\x f_j(\x,\xib) = \nabla_{\x}f_j^o(\x,\w_j^\star,\xib) + \lambda_j^\star\nabla_{\x}f_j^c(\x,\w_j^\star,\xib) \label{gradenv}
	\end{align}
	for all $j = 1, \ldots, m$. If $\ell$ is convex, a convex surrogate of $f$ at a point $\x_0$ is given by $\ell(\f(\x_0,\xib) + \nabla_\x \f(\x_0,\xib)^\top(\x-\x_0),\xib)$ where the $j$-th column of the matrix $\nabla \f(\x,\xib)$ is the gradient in \eqref{gradenv}. Other surrogates are also possible, depending on the structure of $\ell$; see examples in \cite{scutari2016parallel}. We next discuss an application of this idea to precoding in WSNs.

	\subsection{Low complexity precoding in WSNs}\label{secprec}
	\begin{figure}[!h]
		\centering
		\includegraphics[trim=0cm 11cm 22.5cm 1.36cm, clip=true, scale = .7]{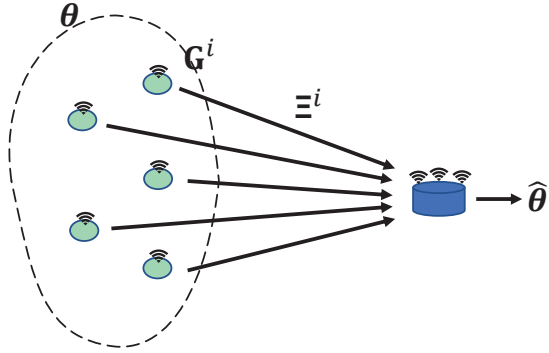} 
		\caption{System model for distributed sensing}
		\label{behfig}
	\end{figure}
	Consider a WSN with $K$ sensors observing a parameter $\so \in \Rn^p$ with zero mean and covariance matrix $\R_\so$. We follow the system model from \cite{behbahani2012linear} depicted in Fig. \ref{behfig}, where the sensors transmit pre-coded versions of received signals in order to improve upon the spectral efficiency. Specifically, the $i$-th sensor observes the $l_i$-dimensional complex signal $\r_i = \H^i\so + \n_i$ where $\H^i \in \Cn^{l_i \times p}$ is the static observation matrix and $\n_i \sim \mathcal{N}(0,\R_{r_i})$ is the independent identically distributed (i.i.d.) noise. Stacking all the sensor observations into an $l = \sum_il_i$-dimensional vector $\r$, we can write the observed signal as $\r = \H\so + \n_r$, where $\H:=[\H^{1\top} ~ \ldots ~ \H^{K\top}]$ and $\r$ and $\n_r$ are stacked versions of $\r_i$ and $\n_{r_i}$, respectively, as in \cite{behbahani2012linear}. The observations at sensor $i$ are linearly precoded with matrix $\G^i \in \Cn^{N_i \times l_i}$ and coherently transmitted to the fusion center (FC) through a coherent MAC (Multiple Acccess Channel) matrix $\Xib^i \in \Cn^{N_F \times N_i}$ where $N_F$ is the number of antennas at the FC. For simplicity, we consider the noiseless case, so that the received signal at the FC is given by $\y_r = \Xib\G\H\so + \Xib\G\n_r$ where $\Xib \in \Cn^{N_F \times \sum_{i}N_i}$ collects the channel gain matrices between the sensors and the FC. Also, here $\G := [\G^1 \oplus \ldots  \oplus \G^K ]$ is a block diagonal $\sum_{i=1}^{K}N_i \times l$ matrix. In this setting, it is well-known from \cite{behbahani2012linear} that $\G$ can be designed so as to ensure that $\hat{\so} = \y_r$, i.e., the parameter estimate is directly recovered from $\y_r$ without the need for equalization. It is remarked that while only a special case of distributed linear precoding is considered here, the proposed approach readily extends to other more general settings. For this example, we let $\xib  = \vec([\mathfrak{Re}(\Xib)~ \mathfrak{Im}(\Xib)])$ so that $d = 2N_F \sum_{i=1}^{K}N_i$ consistent with the notation introduced earlier. 
	
	The mean-square error (MSE) arising within the aforementioned model can be written as \cite{behbahani2012linear}
	\begin{align}
	\zeta(\G,\xib) &= \tr \big( \Xib\G\H\R_{\so}\H^{H}\G^{H}\Xib^{H} - \Xib\G\H\R_{\so} \nonumber \\ 
	& + \Xib\G\R_{n_r}\G^H\Xib^H - \R_{\so}\H^{H}\G^{H}\Xib^{H} + \R_{\so}\big) \label{beforebeh}
	\end{align}
	and the power-constrained MSE minimization problem is given by
	\begin{align}\label{app_prob}
	\min_{\G} ~&~ \zeta(\G,\xib) \\ \text{s. t. } & \hspace{0.35cm} \tr (\G \Ga \G^H) \leq P \nonumber
	\end{align}
	where $P$ is the global power budget and $\Ga = \Ex{\r\r^H} = \H^H\R_{\so}\H+\R_{n_r}$. Alternatively, per-node power constraints of the form $[\G\Ga\G^H]_i \leq P_i$ for $1\leq i \leq K$ may likewise be imposed. Note that since the MSE  depends on the instantaneous channel matrix $\xib$ and therefore the precoder $\G$ must be designed for every coherence interval. Such a requirement may generally be quite challenging, as a closed-form solution to \eqref{app_prob} cannot be found in general. 
	
	Let $\hat{\G}$ denote the precoding matrix corresponding to a specific coherence interval. We consider a realistic scenario where the entries of the channel matrix $\xib$ include both small-scale and large-scale fading components. The large-scale, terrain dependent components are nearly static or vary slowly over time, while the small-scale components contribute to smaller fluctuations occurring at a faster time-scale. Such a split is routinely utilized in the context of wireless communications, especially when modeling the availability of channel state information at various nodes \cite{liu2019robust, zheng2009robust, yang2008distributed, zhang2008statistically}. 
	
	The split allows us to write the precoder matrix  $\hat{\G}$ as $\hat{\G} = \G + \G_\xib$ where $\G$ represents the slowly time-varying component matched to the almost static large-scale fading component while $\G_\xib$ represents the instantaneous component that is different at each coherence interval. If the instantaneous fluctuations in the channel gain are small, we can impose the restriction $\norm{\G_\xib}_F \leq \varepsilon$ for some small $\varepsilon > 0$. The split allows us to express \eqref{app_prob} in the required form of \eqref{ell}-\eqref{fj} by expressing the objective function as
	\begin{align} \label{sec2pro}
	f(\G,\xib) &:= \max_{\G_{\xib}} -\zeta(\G + \G_{\xib},\xib)  \\ 
	\text{s. t. }& \hspace{3.5mm} \varepsilon - \norm{\G_{\xib}}_F \geq 0. \nonumber
	\end{align}
	whose approximate solution can be found easily if $\varepsilon$ is small. It is remarked that such a split is commonly used in the wireless communications literature, where imperfect CSI is often modeled as comprising of the actual channel plus error. It is common to take the CSI error as belonging to a norm ball.
	
	If $\norm{\Ga}_F\leq \omega$, then it follows that the restriction $\tr(\G^H\Ga\G) \leq \tilde{P} :=(P-\omega\varepsilon)^2$ would imply that $\tr((\G+\G_\xib)^H\Ga(\G+\G_\xib)) \leq P$ for all $\xib$. The reduction in the power of $\G$ is required so as to ensure that the overall precoder $\G+\G_\xib$ adheres to the power constraint, regardless of the channel conditions. The optimum value of $\G$ can be found by maximizing $\Ex{f(\G,\xib)}$ under the constraint $\G \in \Xc := \{\tr(\G \Ga \G^H) \leq \tilde{P}\}$. In summary, the per-iteration precoder design problem now reduces to finding $\G_\xib$ for every $\xib$ and subsequently updating $\G$ at every iteration using the SCA or its asynchronous variant proposed in the next section. 
	
	\section{Asynchronous Stochastic SCA} 
	This section details the analysis of the asynchronous stochastic SCA algorithm. The synchronous form of the proposed algorithm is the same as the stochastic successive approximation algorithm in \cite{liu2018stochastic}. The present analysis is however almost entirely different, as it entails characterizing the SCO complexity so as to precisely characterize the effect of delays.
	\subsection{Preliminaries}\label{Asynch_sec}    
	We begin with discussing the required assumptions. 
	\begin{assumption} \label{asreg}
		The set $\text{ri}(\Xc \cap \text{dom}(F) \cap \text{dom}(h))$ is non-empty and $F$ is $L$-smooth. 
	\end{assumption}
	\begin{assumption}\label{asvar}
		The variance of the stochastic gradient is bounded, i.e., $\Ex{\norm{\nabla f(\x,\xib) - \nabla F(\x)}^2} \leq \sigma^2$ for all $\x\in\Xc$. 
	\end{assumption}
	\begin{assumption}\label{assur}
		The surrogate $\fh(\x,\z,\xib)$ of $f(\x,\xib)$ at point $\z$ is $\mu$ strongly convex and $\L$-smooth in $\x$, with $\nabla \fh(\x,\x,\xib) = \nabla f(\x,\xib)$.
	\end{assumption} 
	Assumptions \ref{asreg} and \ref{asvar} are standard in the context of proximal stochastic gradient descent. Assumption \ref{assur} restricts the choice of the surrogate function and was also used in \cite{scutari2016parallel, liu2018stochastic}. Examples of surrogate functions satisfying Assumption \ref{assur} can be found in \cite{scutari2016parallel}. 
	
	In order to characterize the SCO complexity, it is necessary to properly define the $\epsilon$-stationary point of the non-convex problem \eqref{opt-P4}. A point $\bar{\x}$ is a stationary point, if it satisfies $-\nabla F(\bar{\x}) \in \partial (h(\bar{\x}) + \ind_{\Xc}(\bar{\x}))$. A point $\bar{\x}$ is said to be $\epsilon$-stationary if it satisfies:
	\begin{align}
	\min_{\v_{\bar{\x}} \in \partial (h(\bar{\x}) + \ind_{\Xc}(\bar{\x}))} \norm{\v_{\bar{\x}} + \nabla F(\bar{\x})}^2 \leq \epsilon.
	\end{align}
	In general the iterate produced by the algorithm will be random and we will establish their $\epsilon$-stationarity in expectation. 
	
	\subsection{Asynchronous Algorithm}	
	We are now ready to state the proposed asynchronous algorithm. We use $[t]$ to denote a time or iteration such that the delay $t-[t] \in \{0, \ldots, \tau\}$. The proposed algorithm takes the form
	\begin{subequations}\label{algoa}
		\begin{align}
		\x_{t+1} &= (1-\gamma)\x_t + \gamma \xhatb \label{bb}\\
		\y_{t+1} &= (1-\rho)\y_t + \rho\nabla f(\x_t,\xib_t) \label{aa}
		\end{align}
		for positive parameters $\gamma$, $\rho$, and, 
		\begin{align}\label{xhata}
		\xhat &= \arg\min_{\x\in\Xc} h(\x) + \ft(\x,\x_t,\y_t,\xib_t) 
		\end{align}
	\end{subequations}
	where $\ft(\x,\x_t,\y_t,\xib_t):= \rho\fh(\x,\x_t,\xib_t)+(1-\rho)\ip{\y_t,\x}   + (1-\rho)\frac{\mu}{2}\norm{\x-\x_t}^2$ for all $t\geq 1$ and $0 \leq t-\t \leq \tau$. The idea behind these updates is as follows: $\y_{t+1}$ tracks and thus serves as an approximation to the gradient $\nabla F(\x_t)$. The surrogate function $\ft$ is formed from a combination of terms depending on the instantaneous surrogate $\hat{f}$ and the tracked gradient $\y_t$. It can be seen that $\nabla \ft(\x_t,\x_t,\y_t,\xib_t) = (1-\rho)\y_t + \rho\nabla f(\x_t,\xib_t) = \y_{t+1}$ and that $\ft$ is $\mu$ strongly convex in its first argument. Finally, the asynchrony is effected by using $\hat{\x}$ (solution to \eqref{xhata}) from an earlier time $[t]$. This way, the iterates may continue while \eqref{xhata} is being solved, through the use of the old solution. The proposed algorithm is useful when multiple cores or nodes are available for solving \eqref{xhata} in parallel. For instance, if there are $M$ cores available, each can solve \eqref{xhata} over $M$ time slots with $[t] = t-M$ for each $t$, while the overall algorithms continues to run with one update per time slot. The actions to be performed by each core are summarized in Algorithm \ref{algocore} while those being performed by the coordinator are summarized in Algorithm \ref{algofc}. 
	
	In the general case, each core receives some $\x_{t'}$ and $\y_{t'}$ from the central coordinator at time $t'$. Upon receiving these iterates, the core also reads data $\xib_{t'}$ and evaluates $\xh(\x_{t'},y_{t'},\xib_{t'})$ over the next few (at most $\tau$) time slots. Eventually, $\xh(\x_{t'},y_{t'},\xib_{t'})$ is transmitted to the central coordinator at time $t$, which uses it to update $\x_{t+1}$ and $\y_{t+1}$ using \eqref{bb}-\eqref{aa}, and sends them back to the core. The messages received from each core arrive in a first-in-first-out queue, where they wait to be processed. As a result, the updates $(\x_{t+1},\y_{t+1})$ received by each core are unique. 
	
	To ensure consistency, all cores initialize using  $\x_1$ and data point $\xib_1$. The first core to finish calculating $\xh(\x_1,\y_1,\xib_1)$ proceeds to calculate the updated $\x_2$ and starts calculating $\xh(\x_2,\y_2,\xib_2)$. The other cores join in subsequently, calculating $\x_3$, $\x_4, \ldots$ using  $\xh(\x_1,\y_1,\xib_1)$, and rest of the algorithm proceeds as Algorithm \ref{algofc}. In other words, if there are $K$ cores, the initialization ensures that $[t] = 1$ for $1\leq t \leq K$, but the sequence $\{[t]\}_{t=K+1}^T$ has unique entries from $\{ 2, \ldots, T\}$. 
	
	
	
	\begin{algorithm}
		\caption{Operations at each core}
		\label{algocore}
		\textbf{Receive: }$\x$, $\y$ from the coordinator;\\
		\textbf{Observe }$\xib$; \\
		\textbf{Evaluate and send }$\hat{\x}(\x,\y,\xib)$ 
	\end{algorithm}
	
	\begin{algorithm}
		\caption{Asynchronous Stochastic SCA (Asy-SCA) Algorithm} 
		\label{algofc}
		\textbf{Initialize} $\x_1$ , $\y_1$ , $\gamma$ , $\rho$; \\
		\textbf{Transmit }$\x_1$ to all cores;\\
		\For {$t = 1,2,...,T$} {
			\textbf{Receive }$\hat{\x}(\x_{t'},\y_{t'},\xib_{t})$ for some $t'\leq t$;\\
			\textbf{Update }$\x_{t+1}$ using \eqref{bb}; \\
			\textbf{Evaluate }$\nabla f(\x_t,\xib_t)$;\\
			\textbf{Update }$\y_{t+1}$ using \eqref{aa};\\
			\textbf{Transmit }$\x_{t+1}$, $\y_{t+1}$ to the core;  
		}
	\end{algorithm}

	\begin{figure*}[t]
		\centering
		\captionsetup{justification=centering}
		\includegraphics[trim=0cm 3.2cm 0cm 2.1cm, clip=true, width=0.9\linewidth]{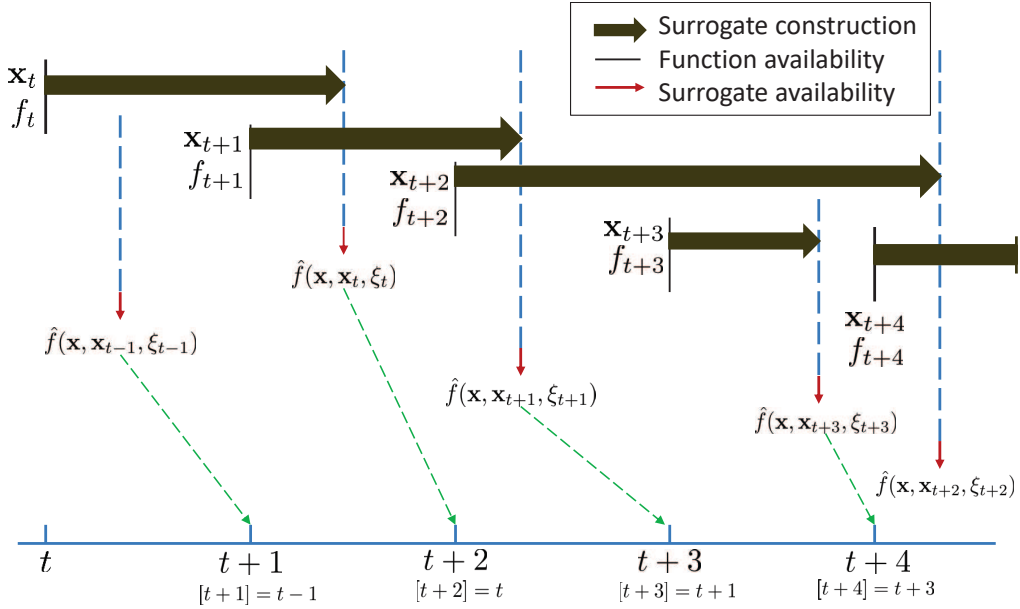}  
		\caption{Execution timeline of the proposed algorithm. Due to surrogate evaluation taking  more than one coherence interval, it is only available at the start of a later slot (depiction inspired from [Fig. 1,  \cite{bedi2019asynchronous}]).}	
		\label{sch}
	\end{figure*} 
	
	Fig. \ref{sch} demonstrates the operation of the AsySCA algorithm. Since the surrogate calculation may not be instantaneous, the surrogate function may only be available at the subsequent time slots, at which point it is obsolete. The proposed asynchronous algorithm is however able to use such obsolete surrogate functions to carry out the updates. For example, the surrogate calculation started at $t + 1$ is completed after $t + 2$, hence surrogate $\hat{f}(\x,\x_{t},\xi_t)$ is available for use at $t + 3$, resulting in a delay of 2 at that node.
	
	\subsection{Complexity analysis} \label{con_pr}
	The SCO complexity analysis for the proposed algorithm is very different from the existing convergence proofs in \cite{liu2018stochastic, yang2016parallel,scutari2016parallel,scutari2016paralel}. For instance, \cite{yang2016parallel} establishes convergence using an asymptotic result from \cite[Lemma 1]{ruszczynski1980feasible}, which cannot be used to obtain the corresponding convergence rates. Instead, the idea here is to develop a recursive relationship on $\norm{\ph_t}^2$ that allows us to construct an appropriate Lyapunov function amenable to telescopic summation. 
	
	Before proceeding, we introduce some compact notation:
	\begin{align}
	\ph_t &:= \y_t - \nabla F(\x_{t-1}) & \Phi_t &:= \Ex{\norm{\ph_t}^2} \label{ph}\\
	\del_t &:= \xhatb - \x_t & \Delta_t &:= \Ex{\norm{\del_t}^2} \label{del}
	\end{align}
	so that the update can be written compactly as $\x_{t+1} = \x_t + \gamma\del_t$. Note that while \eqref{del} holds for $t \geq1 $, \eqref{ph} holds for $t>1$ with  $\ph_1 = \y_1 , \Phi_1 = \Ex{\norm{\y_1}^2} $.
	We begin with determining the precise manner in which $\y_t$ tracks $\nabla F(\x_{t-1})$ in the following lemma. 
	\begin{lemma} \label{lem1}
		For $t \geq 2$, the sequence $\y_t$ satisfies 
		\begin{align}
		\Ek{\norm{\ph_{t+1}}^2} \leq  (1-\rho)\norm{\ph_{t}}^2 + 2\rho^2\sigma^2 + \frac{L^2\gamma^2}{\rho}\norm{\del_{t-1}}^2 \label{trackb}
		\end{align}
		where $\Ek{\cdot}$ is the expectation with respect to $\xib_t$ given $\xib_1, \ldots, \xib_{t-1}$.  Also,	$
		\expect_{\xib_{1}}[\norm{\ph_2}^2] \leq (1-\rho)\norm{\ph_{1}}^2 + 2\rho^2\sigma^2 + \frac{L^2}{\rho}\norm{\nabla F(\x_1)}^2 \label{ph2}
		$.
	\end{lemma}
	Lemma \ref{lem1}, whose full proof is provided in Appendix \ref{appendix:A}, characterizes the single-step evolution of the error incurred in tracking the gradient $\nabla F(\x_t)$.

	\begin{lemma} \label{lem2}
		For $t>1$, given $\eta > 0$, the sequence $U(\x_t)$ satisfies:	
		\begin{align}
		U(\x_{t+1})-U(\x_t) \leq & \frac{\gamma\eta}{2}\norm{\ph_{\t+1}}^2+ \gamma\left(\frac{\gamma L}{2}+\frac{2}{\eta} -\mu\right)\norm{\del_t}^2 \nonumber \\
		&~~~+ \frac{\gamma^3\tau\eta(L^2+\rho^2\L^2+\mu^2)}{2}\sum_{j=1}^{\min \{t,\tau\}} \norm{\del_{t-j}}^2. \label{boundb}
		\end{align} 
	\end{lemma}
	The proof of Lemma \ref{lem2} is provided in Appendix \ref{appendix:B} and utilizes the smoothness of $F$, convexity of $h$ and $\Xc$, and Assumption \ref{assur}. The constant $\eta$ is introduced when applying the Peter-Paul inequality and will be set later. Having established the preliminary results, we are now ready to state the main result of the paper.

	\begin{theorem}  \label{thm1}
		The following bound holds:
		\begin{align}
		\frac{1}{T}\sum_{t=1}^T \Delta_t \leq \frac{\rho\tau\Theta_1 +(1-\rho)\Theta_1 +  2\sigma^2\rho^2T}{\gamma C_{\mu} T}
		\end{align} 
		if $\mu$ is such that $C_{\mu}:= \mu - \frac{\gamma L}{2} - \frac{(1+L^2)\gamma}{\rho} - \tau^2\rho\gamma (L^2+\L^2\rho^2+\mu^2) > 0$ and  $\Theta_1 := \expect[\norm{U(\x_1)}^2] + \expect[\norm{\y_1}^2]$ is a constant.
	\end{theorem}
	
	\begin{IEEEproof}
		Defining $U_t = \Ex{U(\x_t)}$, and taking full expectation in \eqref{trackb} and in the result of Lemma \ref{lem2}, we obtain:
		\begin{align}
		\Phi_{t+1} &\leq  (1-\rho)\Phi_t + 2\rho^2\sigma^2 + \frac{L^2\gamma^2}{\rho}\Delta_{t-1}  \label{trackb2}\\
		U_{t+1}&\leq U_t + \gamma\left(\frac{\gamma L}{2}+\frac{2}{\eta} -\mu\right)\Delta_t +\frac{\gamma\eta}{2}\Phi_{\t+1}  + \frac{\gamma^3\tau\eta C_L}{2}\sum_{j=1}^{\tau} \Delta_{t-j}. \label{boundb2}	
		\end{align}
		where $C_L:=L^2+\L^2\rho^2+\mu^2$. Adding up \eqref{trackb2} and \eqref{boundb2} and defining $\Theta_t := U_t + \Phi_t$, we obtain the recursive relationship:
		\begin{align}
		\Theta_{t+1}&-\Theta_t \leq 2\rho^2\sigma^2 + \frac{L^2\gamma^2}{\rho}\Delta_{t-1} +\frac{\gamma^3\tau\eta C_L}{2}\sum_{j=1}^{\tau} \Delta_{t-j}-\gamma\left(\mu-\tfrac{\gamma L}{2}-\tfrac{2}{\eta}\right)\Delta_t  - \rho\Phi_t + \frac{\gamma\eta}{2}\Phi_{\t+1}.
		\end{align}
		Summing over $t = 1, \ldots, T$, and using the fact that $\Theta_t \geq 0$ for all $t$, we obtain
		\begin{align}
		0 &\leq \rho\tau\Theta_1 + (1-\rho)\Theta_1 + 2\sigma^2 \rho^2 T + \frac{L^2\gamma^2}{\rho}\sum_{t=1}^T\Delta_t +\frac{\gamma^3\tau^2\eta C_L}{2}\sum_{t=1}^T\Delta_t\nonumber\\
		&~~~~~-\gamma\left(\mu-\frac{\gamma L}{2}-\frac{2}{\eta}\right)\sum_{t=1}^T\Delta_t  + (\frac{\gamma\eta}{2}-\rho)\sum_{t=1}^{T+1} \Phi_t \label{tele2}
		\end{align}
		where we have used the inequalities 
		\begin{align}
		\sum_{t=1}^T\sum_{j=1}^{\tau}\Delta_{t-j} &\leq \tau\sum_{t=1}^T\Delta_t \\
		\sum_{t=1}^T \Phi_{\t+1} &\leq \tau\Phi_1 + \sum_{t=2}^{T+1} \Phi_t \leq (\tau-1)\Theta_1 + \sum_{t=1}^{T+1}\Phi_t
		\end{align}
		since $[t] = 1$ for at most $\tau$ time slots and we let $\frac{\gamma\eta}{2}\leq 1$. The last term in \eqref{tele2} can be dropped by choosing $\eta = \frac{2\rho}{\gamma}$. Rearranging, we obtain 
		\begin{align}
		C_{\mu}\frac{1}{T}\sum_{t=1}^T \Delta_t \leq \frac{\rho\tau\Theta_1 +(1-\rho)\Theta_1 +  2\sigma^2\rho^2T}{\gamma T}
		\end{align}
		where $C_{\mu} = \mu - \frac{\gamma L}{2} - \frac{(1+L^2)\gamma}{\rho} - \tau^2\rho\gamma C_L$. In order for this bound to be meaningful, it is required that $\mu$ be such that $C_{\mu} > 0$. 
	\end{IEEEproof}
	
	Theorem \ref{thm1} establishes an important bound on the metric $\frac{1}{T}\sum_{t=1}^T\Delta_t$ and will now be used to establish the required SCO complexity bound. 
	
	\begin{theorem}  \label{the1}
		Let the delay $\tau \leq o(\epsilon^{-1})$, $\mu = \Omega(L)$, and $\L \leq \O(\epsilon^{-1})$. Then, the number of iterations required for at least one iterate $\{\xhat\}_{t\geq \tau+1}$ to be $\epsilon$-stationary on average is $\O\left(\epsilon^{-2}\right)$.
	\end{theorem}
	
	The proof of Theorem \ref{the1} is provided in Appendix \ref{appendix:D}. The SCO complexity bound in Theorem \ref{the1} matches the lower bound for the case when only stochastic gradient information is made available at every iteration $t$ \cite{arjevani2019lower}. Additionally, the delay in Theorem \ref{the1} can be arbitrarily large as long as it is $o(\epsilon^{-1})$. 
	
	It is remarked that for the parameter choices in Theorem \ref{the1}, the classical proximal stochastic gradient descent is not a special case of the proposed algorithm for any choice of the surrogate function. In particular, even for $[t] = t$, neither $\xhat-\x_t$ nor $\y_{t+1}$ are unbiased approximations of $\nabla F(\x_t)$. The proposed algorithm is instead analyzed as a quasi-stochastic gradient algorithm, since $\y_{t+1}$ is still an approximation to $\nabla F(\x_t)$. Note further that the choice $\L \leq \O(\frac{1}{\epsilon})$ allows us to choose a broad class of surrogate functions where $\L \gg \mu$, while still achieving the same worst-case SCO complexity. The maximum allowable delay also need not be bounded by a constant and can be arbitrarily large for $\epsilon$  sufficiently small.

	
	\section{Simulations}\label{sec.sm}
	This section provides simulation results for the low-complexity precoding problem discussed in Sec. \ref{secprec}. We begin with discussing the precoder design approaches considered here and subsequently demonstrate the advantages of the proposed asynchronous algorithm.

	
	
	\subsection{Precoder Design Approaches}\label{precdesign}
	In line with the notation introduced so far, let the coherence intervals be indexed by $t$  and let $\xib_t$ be the channel gains at time $t$.
	
	\subsubsection{Instantaneous Optimal Precoder}
	We begin by stating the optimal approach, where the MSE minimization problem is solved at every coherence interval. That is, the instantaneous precoder is given by
	\begin{align}\label{optimal}
	\hat{\G}_t &:= \arg\min_{\G} \ \ \  \zeta(\G,\xib_t) \\ \text{s. t. }& \hspace{.35cm}\tr (\G \Ga\G^H) \leq P \nonumber 
	\end{align}
	Note that \eqref{optimal} cannot generally be solved in closed form, and solving it numerically also incurs a complexity that is cubic in $K$ and $l$. In practice, such an approach cannot be implemented if the time required to solve \eqref{optimal} is a significant fraction of the coherence time. Nevertheless, we use the optimal solution as a baseline to compare the efficacy of the various suboptimal approaches proposed here. 
	
	
	\subsubsection{Static Optimal Precoder}
	At the other extreme is the static precoder that is optimal in hindsight. Specifically, the  static precoder that minimizes the time-averaged MSE and is given by 
	\begin{align}\label{suboptimal}
	\hat{\G}^s = &\arg\min_{\G} \ \ \ \frac{1}{T}\sum_{t=1}^{T} \zeta(\G,\xib_t)\\
	\text{s. t. }  & \hspace{.35cm} \tr (\G \Ga \G^H) \leq P.  \nonumber
	\end{align}
	Of course, the channel gains $\{\xib_t\}_{t=1}^T$ are not known in advance, so the static precoder is also impractical. In practice however, it may be possible to obtain $\hat{\G}^s$ using an online stochastic approximation algorithm. For instance, if the samples $\xib_t$ are observed sequentially over time, the stochastic gradient descent (SGD) updates take the form $	\hat{\G}^s_{t+1} = \mathcal{P}_{\tr(\G^H\Ga\G)\leq P}\left(\hat{\G}^s_t - \eta_t\nabla \zeta(\G_t,\xib_t)\right)$ where $\mathcal{P}$ denotes the projection operation and $\eta_t$ is the step-size. It is known that given sufficiently large $T$, the iterates $\hat{\G}^s_t$ will approach $\hat{\G}^s$. Note however that even with complete knowledge of $\{\xib_t\}$ in advance, the average MSE of the static precoder will be worse than that of the instantaneous optimal precoder, i.e., $\frac{1}{T}\sum_{t=1}^T\zeta(\hat{\G}^s,\xib_t) > \frac{1}{T}\sum_{t=1}^T \zeta(\hat{\G}_t,\xib_t)$ since the static precoder is not adapted to the channel variations. Moreover, the projection operation in SGD is not much simpler than solving the instantaneous problem \eqref{optimal}. The precoder designs proposed next will remedy this weakness and therefore outperform the static precoder. 
	
	\subsubsection{Hybrid Precoder via Envelope Theorem}
	When the channel variations are sufficiently small, it may be possible to split the precoder as $\hat{\G} = \G + \G_\xib$. As explained in \ref{secprec}, the static part of the hybrid precoder is given by the solution to:
	\begin{align}\label{p1}
	\G^h &= \arg\max_{\G} \expect[f(\G,\xib)] \tag{$\mathcal{P}_1$} \\
	\text{s. t. }& \hspace{3.5mm} \tr(\G^H\Ga\G) \leq \tilde{P} \nonumber
	\end{align}
	where the objective function is as specified in \eqref{sec2pro}. The proposed algorithm is applied on \ref{p1} and the iterate $t$ of $\G$ is represented by $\G_t$. So, at time $t$, given the current iterate $\G_t$, we first solve the following per coherence interval problem:
	\begin{align} \label{g1}
	\G_{\xib_t} &= \arg\max_{\G} -\zeta (\G_t + \G,\xib_t) \\
	\text{s. t. } & \hspace{5mm} \norm{\G}_F \leq \varepsilon \nonumber
	\end{align}
	If $\varepsilon$ is small, one option is to approximate $\zeta(\G_t+\G_{\xib_t},\xib_t) \approx \zeta(\G_t,\xib_t) + \ip{\nabla\zeta(\G_t,\xib_t),\G_{\xib_t}}$ so that 
	\begin{align}\label{approxG}
	\G_{\xib_t} \approx -\frac{\nabla\zeta(\G_t,\xib)}{\norm{\nabla\zeta(\G_t,\xib)}_F}\varepsilon.
	\end{align}
	Alternatively, other approximations exploiting the structure of the MSE function may be considered. Finally, the surrogate for $f$ is given by 
	\begin{align}\label{fhatG}
	\hat{f}(\G,\G_t,\xib_t) = &-\zeta(\G_t,\xib_t) - \ip{\nabla\zeta(\G_t+\G_{\xib_t},\xib_t),\G-\G_t} \nonumber\\
	&~~~~~+ \frac{\mu}{2}\norm{\G-\G_t}_F^2.
	\end{align}
	Application of the proposed Asy-SCA updates will therefore allow us to obtain $\G_{t+1}$. 
	
	The proposed algorithm overcomes the limitations of both the instantaneous and static optimal precoder designs. Specifically, all nodes maintain and update $\G_t$ at every coherence interval. At time $t$, we solve \eqref{g1} approximately, and use $\G_t+\G_{\xib_t}$ for transmission. However, the construction of the surrogate function as well as the update of $\G_t$ (which requires solving an optimization problem) need not occur immediately, but may instead be carried out whenever the required computations have been completed, thanks to the asynchronous updates allowed by the proposed algorithm.

	\subsubsection{Hybrid Precoder for Convex Case}
	We detail yet another hybrid approach that is applicable if the MSE is a convex quadratic function of $\G$, allowing us to construct surrogates without the use of the envelope theorem in Sec. \ref{secprec}. We begin with observing that if $\varepsilon$ is small, the MSE expression can be approximated as $\zeta(\G+\G_{\xib},\xib) \approx \zeta(\G,\xib) + \ip{\nabla\zeta(\G,\xib),\G_{\xib}}$. Therefore, the optimal value of $\G_{\xib}$ that minimizes the MSE for a given value of $\G$ can be approximated as in \eqref{approxG}. Now, instead of using the envelope theorem in \eqref{fhatG}, we construct a surrogate by making use of the convex quadratic form of the MSE. 
	
	For notational brevity, let $\g$ and $\g_\xib$ denote the vectorized versions of $\G$ and $\G_{\xib}$, and let the MSE expression be given by $\zeta(\g,\xib) = \g^H\A_{\xib}\g - \b_{\xib}^H\g - \g^H\b_{\xib}^H + q_{\xib}$ where $\A_\xib$ is positive semi-definite. Therefore, the MSE approximation is given by $\zeta(\g+\g_\xib,\xib) \approx \zeta(\g,\xib) + \d^H(\g,\xib)\g_{\xib}$ where  $\d(\x,\xib):= \nabla \zeta(\g,\xib) =2(\A_\xib\g - \b_\xib)$. As in \eqref{approxG}, the $\g_{\xib}$ that minimizes the MSE at time $t$ would be given by $	\g_\xib = -\frac{\d(\g,\xib)}{\norm{\d(\g,\xib)}_F}\varepsilon$, which, upon substituting into the original MSE expression, yields
	\begin{align}
	\zeta(\g,\xib) &= \g^H\A_\xib\g - \b^H_\xib\g -\g^H\b_{\xib}^H + q_\xib - \varepsilon\norm{\d(\g,\xib)} \nonumber \\
	& + \varepsilon^2\frac{\d^H(\g,\xib)\A\d(\g,\xib)}{\norm{\d(\g,\xib)}^2}.
	\end{align}
	Here observe that if $\A$ is positive semi-definite, the first four terms constitute the convex component of the objective while the last two terms could be nonconvex, i.e., $\zeta(\g,\xib) = \zeta^c(\g,\xib) + \zeta^{\bar{c}}(\g,\xib)$, where
	\begin{subequations}
		\begin{align}
		\zeta^c(\g,\xib) &:= \g^H\A_\xib\g - \b_\xib^H\g - \g^H\b_{\xib}^H + q_\xib \label{con} \\
		\zeta^{\bar{c}}(\g,\xib) &:= \frac{\ep^2\d^H(\g,\xib)\A_\xib\d(\g,\xib)}{\norm{\d(\g,\xib)}^2} - \ep\norm{ \d(\g,\xib)} \label{noncon} 
		\end{align}
	\end{subequations}
	Finally, the surrogate function at $\g_t$ can be constructed as $\zeta^c(\g,\xib_t) + \zeta^{\bar{c}}(\g_t,\xib_t) + \ip{\nabla \zeta^{\bar{c}}(\g_t,\xib_t),\g-\g_t} + \tfrac{\mu}{2}\norm{\g-\g_t}^2$. It is remarked that the function $\norm{\d(\g,\xib)}$ is not differentiable, but can be smoothened by replacing it with $\sqrt{\norm{\d(\g,\xib)}^2+\upsilon^2}-\upsilon$ for some $0 < \upsilon \ll 1$; see \cite{beck2017first} (Example 10.44). The results obtained in the next section will use this approximation. 
	
	Summarizing, the two hybrid precoders considered here circumvent the issues faced by classical precoders. In each of the two precoders, the instantaneous component $\G_\xib$ must still be evaluated at every time instant, but the same is possible through the use of approximations. The static component, on the other hand, can be updated in an asynchronous fashion, thereby tolerating delays arising due to computational issues. 
	
	\subsection{Numerical Results} \label{Results}
	We begin with detailing the simulation setting considered here. For simplicity, we consider a small network $p = K = N_{FC} = N_{i} = l_{i} = 2$ for $1\leq i \leq K$, though the proposed algorithm applies to larger networks provided that $p < N = \sum_{i=1}^K N_i$. Both real and imaginary entries of the observation matrix $\H$ are randomly selected to be $+1$ or $-1$, while the covariance $\R_{\so}$ of the parameter vector and that of the sensor noise $\R_{\n_r}$ are taken to be proportional to identity. The sensor noise is assumed to be 30 dB below the received signal power. We select the total power budget $P = 10$ units. 
	
	To simulate channels with low variability, for each Monte-Carlo run, we first select a matrix $\Xib_0$ with zero mean and complex Gaussian distributed entries. At subsequent coherence intervals or time slots, we simply select channel matrices by adding complex Gaussian random variables to $\Xib_0$ with a pre-set standard deviation. Such a process allows us to precisely control the channel variations and study their impact on the algorithm performance. Unless otherwise specified, we assume that the channel standard deviation of the norm of $\Xib$ is  $0.05$, while the channel norm is unity on average. All results are averaged over 200 Monte-Carlo runs. 
	
	In general, initialization is critical for iterative non-convex algorithms. In order to ensure proper initialization, we calculate the static precoder using five random channel matrices. That is, the initial precoder $\G_0 := \arg\min_{\G} \sum_{i=1}^5 \zeta(\G,\xib'_i)$ subject to $\tr(\G\Ga\G^H) \leq P$, where $\{\xib'_i\}$ are arbitrary channel observations. In practice and for a fair comparison, we assume that the proposed algorithm ``starts'' at the $6$-th time slot. Finally, for each of the two hybrid algorithms, there are four parameter that must be tuned, namely, $\varepsilon$, $\mu$, $\rho$, and $\gamma$. These parameters were hand-tuned on a separate channel matrix and subsequently retained for all the experiments. The proposed algorithms were relatively robust to the choices of $\mu$, $\rho$, and $\gamma$, but the performance did depend on the choice of $\varepsilon$, which encodes the amount of channel variations. For the hybrid precoder designed using envelope theorem, we used  $\ep =0.05,    \mu = 0.2,	\rho = 0.1$, and $\gamma = 0.001$. For the hybrid precoder designed for the convex case, 	$\ep = 0.02$, $\mu = 0.01$, $\rho = 0.01$, and $\gamma = 0.001$. 
	
	We consider a setting where the minimization of the surrogate takes between up to five coherence intervals, translating to a maximum delay of $\tau = 5$ time slots. The time taken to carry out the update is random. 
	
\begin{figure}[h]		
		\centering
		\captionsetup{justification=centering}
		\includegraphics[width=0.75\linewidth]{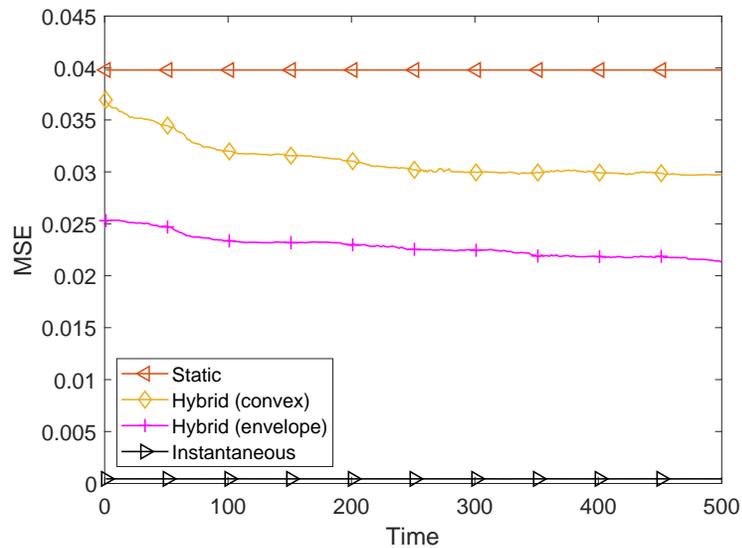}
		\caption{MSE vs time ($t$) using Asy-SCA  for the various precoder design approaches}
		\label{conver}		
	\end{figure}

	We begin with studying the performance of the proposed precoder design approach and contrasting it with that of the instantaneous and static precoders. Fig. \ref{conver} shows the evolution of the MSE of the four schemes outlined in Sec. \ref{precdesign}. First observe that as expected, the instantaneous precoder incurs the lowest MSE, since it aligns its precoder exactly to the exact channel observed at every coherence interval. Interestingly however, the proposed Asy-SCA algorithms outperform the static precoder, even though the static precoder is designed using full knowledge of the future channel matrices. The hybrid precoder design based on the envelope theorem performs significantly better than the static precoder, while running in an online fashion and handling computational delays.

	Next, we compare the proposed Asy-SCA algorithm with its synchronous version CS-SCA, which is known to converge asymptotically from \cite{liu2018stochastic} and \cite{liu2018online}. For a proper comparison, we consider two variants of CS-SCA algorithm. 
	\begin{itemize}
		\item The \textbf{genie-aided CS-SCA} algorithm, where each update takes exactly one coherence interval. As explained earlier, this essentially means that the minimization of the surrogate function as well as the other updates occur almost instantaneously so that the results can be used for precoding in the current coherence interval itself. The genie-aided CS-SCA is also equivalent to running the proposed algorithms $\tau = 0$. 	
		\item The \textbf{practical CS-SCA} algorithm, where surrogate minimization takes between zero and five coherence intervals. But since the CS-SCA algorithm cannot use stale gradients, it must wait till the surrogate is minimized. Note that in this case, the updates will lag behind $t$ and the memory requirements increase linearly with $T$.
	\end{itemize}
	
	\begin{figure}[!h]
		
		\centering
		\captionsetup{justification=centering}
		\includegraphics[width=.75\linewidth]{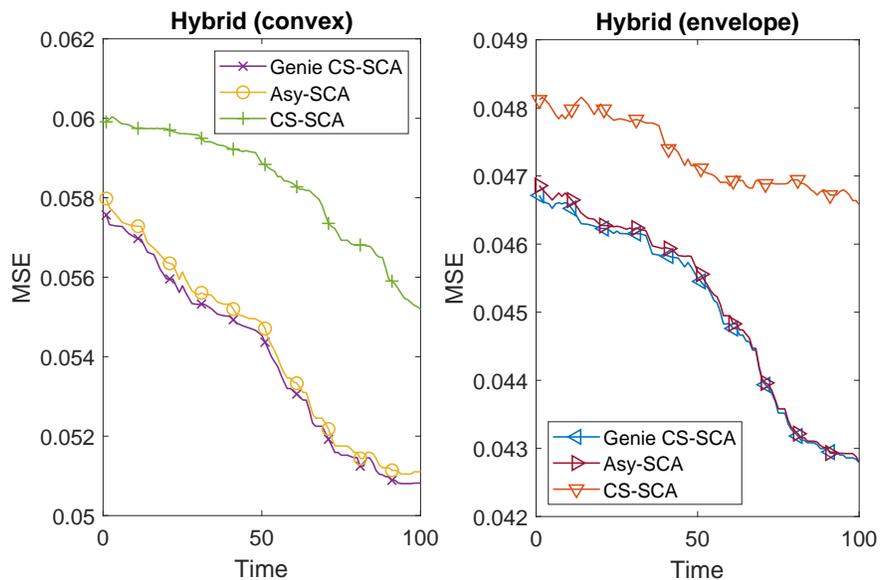}
		\caption{MSE comparison of the genie and practical variants of CS-SCA \cite{liu2018stochastic}, as well as the proposed Asy-SCA algorithm with $\tau = 5$. }
		\label{comp1}
	\end{figure}
	Fig. \ref{comp1} shows the evolution of the MSE of the proposed and the two CS-SCA variants. It can be seen that the proposed Asy-SCA algorithm is quite close to that of the genie-aided CS-SCA even though it allows delayed updates. Indeed, the practical CS-SCA, where the delays result in nodes waiting to carry out the updates lags behind considerably, since it carries out fewer updates per time slot. Fig. \ref{comp1} demonstrates the gains that can be obtained by allowing delayed SCA updates. 
	
\begin{figure}[h]		
		\centering
		\captionsetup{justification=centering}
		\includegraphics[width=0.75\linewidth]{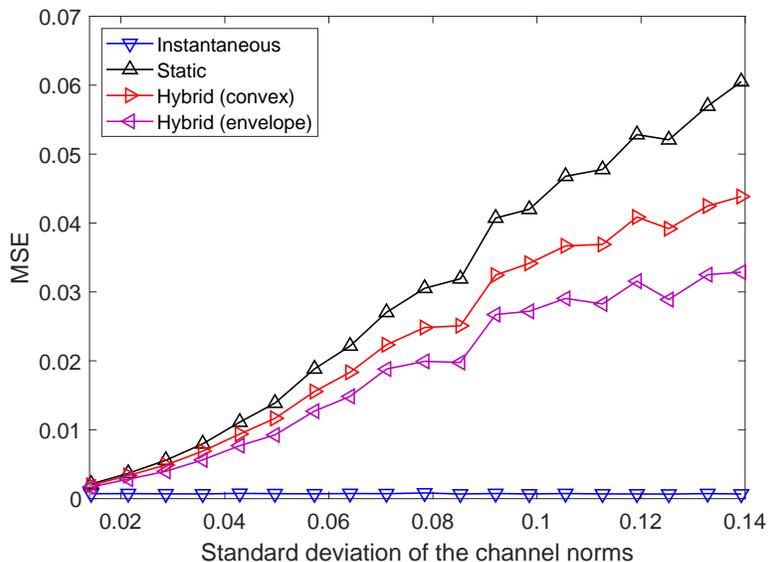}
		\caption{Plot of MSE (seen at 500-th time slot) with standard deviation of the channel norms.}
		\label{comp2}		
	\end{figure}

	It is instructive to study the effect of channel variations, since the proposed approaches utilized various approximations to speed-up the computations. Fig. \ref{comp2} shows the plot of steady-state MSE achieved by the different approaches for different values of the standard deviation. As expected the MSE of the instantaneous precoder remains the same irrespective of the channel variations. However, the MSE of the other three approaches increase. However, MSE of the proposed approach increases to a lesser extent, and continues to give improvement over static precoder design even for approximately 15$\%$ channel variations. 
	
	\section{Conclusion}
	This work studies the non-asymptotic performance of the stochastic successive convex approximation (SCA) algorithm for problems with non-convex loss functions and convex constraint and regularizers. The analysis reveals that the stochastic SCA algorithm has the same iteration complexity as the proximal SGD. The insight allows us to develop an asynchronous variant of the stochastic SCA algorithm where surrogates calculated at earlier iterations can be used. We characterize a bound on the maximum delays that algorithm can tolerate while still achieving the same order of the iteration complexity. Subsequently, the proposed algorithm is applied to the problem of designing precoding matrices in wireless sensor networks utilizing a coherent medium access control mechanism. While the precoder design is essentially a per-iteration problem, we assume that the channel is nearly static, so that the time-varying component of the precoding matrix is small and can be approximated using a closed-form solution. The static component on the other hand is learned in an online fashion using the proposed asynchronous stochastic SCA algorithm. Detailed simulation results demonstrate the efficacy of the proposed framework.

	
	\appendices
	
	\section{Proof of Lemma \ref{lem1}}\label{appendix:A}	
	The proof of Lemma \ref{lem1} follows along the lines of \cite{wang2017stochastic} (Lemma 1). Define $\w_t := (1-\rho)(\nabla F(\x_t) - \nabla F(\x_{t-1}))$, and note that since $\nabla F$ is $L$-Lipschitz, we have that $\norm{\w_t} \leq L(1-\rho)\norm{\x_t-\x_{t-1}}$. Therefore, we have
	\begin{align}
	\ph_{t+1} + \w_t &= (1-\rho)\y_t + \rho\nabla f(\x_t,\xib_t) - \nabla F(\x_t) + (1-\rho)(\nabla F(\x_t) - \nabla F(\x_{t-1}))\\
	&=(1-\rho)\ph_{t} + \rho(\nabla f(\x_t,\xib_t) - \nabla F(\x_t))
	\end{align}
	Taking squared norm and expectation with respect $\xib_t$, we obtain
	\begin{align}
	&\Ek{\norm{\ph_{t+1} + \w_t}^2} = \Ek{\norm{(1-\rho)\ph_{t} + \rho(\nabla f(\x_t,\xib_t) - \nabla F(\x_t))}^2}. 
	\end{align}
	Expanding the right-hand side, observe that the cross-term vanishes since $
	\Ek{\nabla f(\x_t,\xib_t) - \nabla F(\x_t)}  = 0$. Therefore, we have
	\begin{align}
	&\Ek{\norm{\ph_{t+1} + \w_t}^2} = (1-\rho)^2\norm{\ph_{t}}^2 + \rho^2\Ek{\norm{\nabla f(\x_t,\xib_t) - \nabla F(\x_t)}^2} 
	\end{align}
	Using Peter-Paul inequality, we have that
	\begin{align}
	\Ek{\norm{\ph_{t+1}}^2} &\leq (1+\rho)\Ek{\norm{\ph_{t+1} + \w_t}^2} + \left(1+\tfrac{1}{\rho}\right)\norm{\w_t}^2\label{ppaul}\\
	&\leq (1+\rho)(1-\rho)^2\norm{\ph_{t}}^2 + (1+\rho)\rho^2\Ek{\norm{\nabla f(\x_t,\xib_t) - \nabla F(\x_t))}^2} \nonumber\\
	&~~~~~~~~+ L^2\frac{1+\rho}{\rho}(1-\rho)^2\norm{\x_t-\x_{t-1}}^2\\
	&\leq (1-\rho)\norm{\ph_{t}}^2 + 2\rho^2\sigma^2 + \frac{L^2}{\rho}\norm{\x_t-\x_{t-1}}^2 \\
	&= (1-\rho)\norm{\ph_{t}}^2 + 2\rho^2\sigma^2 + \frac{L^2\gamma^2}{\rho}\norm{\del_{t-1}}^2
	\end{align}
	where we have used the fact that $\x_{t} = \x_{t-1} + \gamma\del_{t-1}$. Similarly the bound on $\expect_{\xib_{1}}[\norm{\ph_2}^2]$ can be obtained using \eqref{ppaul} and $\w_1 = (1-\rho)\nabla F(\x_1)$.
	
	\section{Proof of Lemma \ref{lem2}}\label{appendix:B}
	We first state two key inequalities that will allow us to establish Lemma \ref{lem2}. 
	\begin{lemma}
		It holds that
		\begin{subequations}
			\begin{align}
			&-\ip{(1-\rho)\y_\t + \v_\t + \rho \nabla \fh(\x_t,\x_\t,\xib_\t), \del_t}\geq \mu\norm{\del_t}^2 + (1-\rho)\mu\ip{\x_t - \x_\t,\del_t}. \label{funda}\\
			&F(\x_{t+1})	\leq F(\x_t) + h(\x_t) - h(\x_{t+1}) - \gamma\ip{\ph_{\t+1}, \del_t} + \gamma\rho \ip{\nabla f(\x_\t,\xib_\t) - \nabla \fh(\x_t,\x_\t,\xib_\t),\del_t} \nonumber\\
			& ~~~~~~~~~~~~~~+ \gamma\ip{\nabla F(\x_t) - \nabla F(\x_\t),\del_t}+  \gamma\left(\tfrac{\gamma L}{2}-\mu\right)\norm{\del_t}^2 + \gamma\mu(1-\rho) \ip{\x_\t - \x_t,\del_t}.\label{pp0}
			\end{align}
		\end{subequations}
	\end{lemma}
	
	\begin{IEEEproof}[Proof of \eqref{funda}]
		Denoting $\v_\t \in \partial h(\xhatb)$, the optimality condition for \eqref{xhata} at time $[t]$ can be written as
		\begin{align}
		(1-\rho)&\ip{\y_\t + \mu(\xhatb - \x_\t),\x - \xhatb} \nonumber\\
		&~~~~+ \rho\ip{\nabla \fh(\xhatb,\x_\t,\xib_\t),\x - \xhatb}  \nonumber\\
		&~~~~~~~~~+\ip{\v_\t,\x - \xhatb} \geq 0 
		\end{align}
		for any $\x \in \Xc$. Replacing $\x$ with $\x_t$, and rearranging, we obtain
		\begin{align}
		&-\ip{(1-\rho)\y_\t + \v_\t+\rho\nabla \fh(\xhatb,\x_\t,\xib_\t),\del_t} \nonumber 
		\\
		&~~~~~~~~~~~~~~~~~~~~~~~~~~~~~~~~~\geq (1-\rho)\mu\ip{\xhatb - \x_\t,\del_t}.
		\end{align}
		where recall that $\del_t= \xhatb-\x_t$. Adding and subtracting $\rho \nabla \fh(\x_t,\x_\t,\xib_\t)$ inside the left-hand side, we have that
		\begin{align}
		-\ip{(1-\rho)\y_\t &+ \v_\t + \rho \nabla \fh(\x_t,\x_\t,\xib_\t), \del_t} \nonumber\\
		&~~~~  + \rho\ip{\nabla \fh(\x_t,\x_\t,\xib_\t)-\nabla \fh(\xhatb,\x_\t,\xib_\t),\del_t} \nonumber\\
		&~~~~~~~~~~\geq (1-\rho)\mu\ip{\del_t + \x_t - \x_\t,\del_t}. 
		\end{align}
		Since $\fh(\cdot,\x_\t,\xib_\t)$ is $\mu$-strongly convex, the strong monotonicity of $\nabla \fh(\cdot, \x_\t, \xib_\t)$ implies that 
		\begin{align}
		\ip{\nabla \fh(\xhatb,\x_\t,\xib_\t) - \nabla &\fh(\x_t,\x_\t,\xib_\t),\del_t}\geq \mu \norm{\del_t}^2
		\end{align}
		which upon substituting yields
		\begin{align}
		&-\ip{(1-\rho)\y_\t + \v_\t + \rho \nabla \fh(\x_t,\x_\t,\xib_\t), \del_t}\geq \mu\norm{\del_t}^2 + (1-\rho)\mu\ip{\x_t - \x_\t,\del_t}.
		\end{align}
		which is the required inequality. 		
	\end{IEEEproof}
	
	\begin{IEEEproof}[Proof of \eqref{pp0}]
		Since $h$ is convex, we have that
		\begin{align}
		h(\x_{t+1}) &\leq \gamma h(\xhatb) + (1-\gamma)h(\x_t) \\
		&\leq \gamma(h(\x_t) + \ip{\v_\t,\del_t}) + (1-\gamma)h(\x_t)\\
		\Rightarrow \ip{\v_\t,\del_t} &\geq \frac{1}{\gamma}(h(\x_{t+1}) - h(\x_t)) . \label{hdiff}
		\end{align}
		Substituting $\y_{\t+1} = (1-\rho)\y_\t + \rho \nabla f(\x_\t,\xib_\t)$ into \eqref{funda} and using \eqref{hdiff}, we obtain
		\begin{align}
		&\ip{-\y_{\t+1} - \rho \nabla \fh(\x_t,\x_\t,\xib_\t) + \rho \nabla f(\x_\t,\xib_\t), \del_t} \\
		&~~~~~~~\geq \mu\norm{\del_t}^2 + \mu(1-\rho)\ip{\x_t - \x_\t,\del_t} +  \tfrac{1}{\gamma}(h(\x_{t+1}) - h(\x_t)).\nonumber
		\end{align}
		Adding and subtracting $\nabla F(\x_\t) - \nabla F(\x_t)$, we obtain
		\begin{align}
		\ip{\nabla F(\x_\t) - \y_{\t+1}, \del_t} &- \ip{\nabla F(\x_t),\del_t}  + \ip{\nabla F(\x_t) - \nabla F(\x_\t),\del_t}\nonumber\\
		&+ \rho \ip{\nabla f(\x_\t,\xib_\t) - \nabla \fh(\x_t,\x_\t,\xib_\t),\del_t} \nonumber\\
		&~~~~~~~~~\geq  \mu\norm{\del_t}^2 + \mu(1-\rho)\ip{\x_t - \x_\t,\del_t} +  \tfrac{1}{\gamma}(h(\x_{t+1}) - h(\x_t)).
		\end{align}
		which, upon rearranging, yields
		\begin{align}
		\gamma\ip{\nabla F(\x_t),\del_t} &\leq h(\x_t) - h(\x_{t+1}) - \mu\gamma\norm{\del_t}^2 - \mu\gamma(1-\rho)\ip{\x_t - \x_\t,\del_t} \nonumber\\
		&\hspace{-2cm}+ \gamma\rho \ip{\nabla f(\x_\t,\xib_\t) - \nabla \fh(\x_t,\x_\t,\xib_\t),\del_t} + \gamma\ip{\nabla F(\x_t) - \nabla F(\x_\t),\del_t} - \gamma\ip{\ph_{\t+1}, \del_t}. \label{gradF}
		\end{align}
		
		To obtain the required inequality, we use the $L$-smoothness of $F$ which yields
		\begin{align}
		F(\x_{t+1}) & \leq F(\x_t) + \ip{\nabla F(\x_t), \x_{t+1}-\x_t} + \frac{L}{2}\norm{\x_{t+1}-\x_t}^2\nonumber\\
		&\leq F(\x_t) + \gamma\ip{\nabla F(\x_t), \del_t} + \frac{\gamma^2L}{2}\norm{\del_t}^2 \label{gradE}
		\end{align}
		Substituting \eqref{gradF} into \eqref{gradE}, we obtain \eqref{pp0}. 
	\end{IEEEproof}
	
	Having established the basic results, we utilize the Peter-Paul inequality with parameter $\eta$ on all the terms in \eqref{pp0} involving an inner product with $\del_t$:
	\begin{align}
	\ip{\ph_{\t+1}, \del_t} &\leq \frac{1}{2\eta}\norm{\del_t}^2 + \frac{\eta}{2}\norm{\ph_{\t+1}}^2 \label{pp1}\\
	\rho \ip{\nabla \fh(\x_\t,\x_\t,\xib_\t) - \nabla \fh(\x_t,\x_\t,\xib_\t),\del_t} &\leq \frac{1}{2\eta}\norm{\del_t}^2 + \frac{\rho^2\eta \L^2}{2}\norm{\x_t-\x_\t}^2 \label{pp2}\\
	\ip{\nabla F(\x_t) - \nabla F(\x_\t),\del_t} &\leq \frac{1}{2\eta}\norm{\del_t}^2 + \frac{\eta L^2}{2}\norm{\x_t-\x_\t}^2 \label{pp3}\\
	\mu(1-\rho) \ip{\x_\t - \x_t,\del_t} &\leq  \frac{1}{2\eta}\norm{\del_t}^2 + \frac{\eta \mu^2(1-\rho)^2}{2}\norm{\x_t-\x_\t}^2 \label{pp4a}
	\end{align}
	where we have also used the $\L$-smoothness of $\fh$ and $L$-smoothness of $F$. 	
	
	Further, since $t-\t\leq \tau$, $\x_{t+1}-\x_t = \gamma\del_t$, it holds that
	\begin{align}
	\norm{\x_t-\x_\t}^2 &= \norm{\sum_{j=0}^{t-\t-1}(\x_{t-j} - \x_{t-j-1})}^2\leq \tau\sum_{j=0}^{\tau-1}\norm{\x_{t-j}-\x_{t-j-1}}^2  \leq \gamma^2\tau\sum_{j=1}^{\tau} \norm{\del_{t-j}}^2
	\end{align}
	where the summation only includes terms for which $j\leq t-1$. Substituting into \eqref{pp4a} and using the fact that $\rho \geq 0$,  we obtain
	\begin{align}
	\mu(1-\rho)& \ip{\x_\t - \x_t,\del_t} \leq \frac{1}{2\eta}\norm{\del_t}^2 + \frac{\eta \mu^2\gamma^2\tau}{2}\sum_{j=1}^{\tau} \norm{\del_{t-j}}^2 \label{pp4}
	\end{align}
	Therefore, substituting \eqref{pp1}-\eqref{pp4} into \eqref{pp0}, we obtain the required result.

	\section{Proof of Theorem \ref{the1}}\label{appendix:D} 
	Recall that $\del_t := \xhatb - \x_t$. For this proof, we also define $\dt := \xhatb-\x_\t$ which can be bounded as
	\begin{align}
	\norm{\dt}^2 \leq 2\norm{\del_t}^2 + 2\gamma^2\tau\sum_{j=1}^\tau \norm{\del_{t-j}}^2.
	\end{align} 
	From the update equation, we have that for all $t \geq 1$, there exists $\vc_\t \in \Hc_\t := \partial (h + \ind_{\Xc})\mid_{\x = \xhatb}$ such that
	\begin{align}
	&(1-\rho)\y_\t + (1-\rho)\mu\dt + \vc_\t+ \rho\nabla \fh(\xhatb,\x_\t,\xib_\t) = 0. 
	\end{align}
	In order to characterize the SCO complexity, let us also define $\Pi_t:=\Ex{\min_{\v \in \Hc_t}\norm{\v  + \nabla F(\xhat)}^2}$.
	
	Since a data point $\xib_t$ is used at most once over $T$ iterations, we have that
	\begin{align}
	\sum_{t=\tau+1}^T\Pi_t&=\sum_{t=\tau+1}^T\Ex{\min_{\v \in \Hc_t}\norm{\v  + \nabla F(\xhat)}^2} \nonumber\\
	&\leq \sum_{t=\tau+1}^T\Ex{\norm{\vc_t  + \nabla F(\xhat)}^2}\nonumber\\
	&\leq \sum_{t=1}^T\Ex{\norm{\vc_\t  + \nabla F(\xhatb)}^2} 
	\end{align}	
	Therefore, using the update equation and substituting $\y_{\t+1} = (1-\rho)\y_\t + \rho\nabla \fh(\x_\t,\x_\t,\xib_\t)$, we can write:	
	\begin{align}
	&\vc_\t  + \nabla F(\xhatb) = \nabla F(\xhatb) \nonumber \\
	& - \nabla F(\x_\t)  + \nabla F(\x_\t) - \y_{\t+1} -(1-\rho)\mu\dt \label{cor1}\\
	&+ \rho (\nabla \fh(\x_\t,\x_\t,\xib_\t) - \nabla \fh(\xhatb,\x_\t,\xib_\t)) \nonumber
	\end{align}
	The term-differences in \eqref{cor1} can be bounded as
	\begin{align}
	\nabla F(\xhatb) - \nabla F(\x_\t) &\leq L\norm{\dt} \\
	\nabla F(\x_\t) - \y_{\t+1} & \leq \norm{\ph_{\t+1}} \\
	\nabla \fh(\xhatb,\x_\t,\xib_\t) - \nabla \fh(\x_\t,\x_\t,\xib_\t) &\leq \L\norm{\dt}
	\end{align}
	Substituting these bounds, taking expectation, and using norm inequalities, we obtain
	\begin{align}
	\frac{1}{T-\tau}\sum_{t=\tau+1}^T\Pi_t &\leq \frac{2(L+\L\rho + \mu(1-\rho))^2}{T-\tau}\sum_{t=1}^T\Ex{\norm{\dt}^2} +\frac{2}{T-\tau}\sum_{t=1}^{T-\tau}\Phi_{\t+1} 
	\end{align}
	which implies that
	\begin{align}
	\min_{\tau+1\leq t \leq T}\Pi_t &\leq \frac{4(1+\gamma^2\tau)}{T-\tau}(L+\L\rho + \mu(1-\rho))^2\sum_{t=1}^{T-\tau}\Delta_t + \frac{2}{T-\tau}\sum_{t=1}^T \Phi_t \label{almt}
	\end{align}
	The second term can be bounded by taking expectation in \eqref{trackb}, and summing over $t = 1, \ldots, T$, we obtain
	\begin{align}
	\Phi_{T+1}\leq \Phi_1 - \rho\sum_{t=1}^T \Phi_t + 2\rho^2\sigma^2T + \frac{L^2\gamma^2}{\rho}\sum_{t=1}^T \Delta_t
	\end{align}
	which yields
	\begin{align}
	\sum_{t=1}^T \Phi_t \leq \frac{\Phi_1}{\rho} + 2\rho\sigma^2T + \frac{L^2\gamma^2}{\rho^2}\sum_{t=1}^T \Delta_t  \label{phibound}
	\end{align}
	Substituting \eqref{phibound} into \eqref{almt} and using $T-\tau \geq T/2$, we obtain
	\begin{align}
	\min_{\tau+1\leq t \leq T}\Pi_t &\leq 8(1+\gamma^2\tau)(\frac{L^2\gamma^2}{\rho^2}+(L+\L\rho + \mu(1-\rho))^2)\frac{1}{T}\sum_{t=1}^T\Delta_t +\frac{4\Phi_1}{\rho T} + 8\rho\sigma^2 \\
	& \leq 8(1+\gamma^2\tau)\left ( \frac{L^2\gamma^2}{\rho^2}+ (L+\L\rho + \mu(1-\rho))^2\right) \frac{\rho\tau\Theta_1+ (1-\rho)\Theta_1 +2\sigma^2\rho^2T  }{C_{\mu}T\gamma}
	\nonumber\\
	&~~~~~~~~~~+\frac{4\Phi_1}{\rho T} + 8\rho\sigma^2 \label{longform}
	\end{align}
	
	Let us choose $\gamma = \rho = T^{-\frac{1}{2}}$. Suppose further that $\tau = o(T^{\frac{1}{2}})$. Then, if we can ensure that $\mu = L$ and $\L \leq \O(\sqrt{T})$, we would have that $C_\mu = L$, and keeping only the dominant term in \eqref{longform} yields:
	\begin{align}
	\min_{\tau+1\leq t \leq T}\Pi_t \leq \O\left(\frac{1}{\sqrt{T}}\right)
	\end{align}
	Equivalently, the SCO complexity is given by $\O(\epsilon^{-2})$.

	\ifCLASSOPTIONcaptionsoff
	\newpage
	\fi
	
	\bibliographystyle{IEEEtran}
	\bibliography{IEEEabrv,citation}

\end{document}